\documentclass[10pt]{article}

\usepackage[authoryear,round]{natbib}

\usepackage[utf8]{inputenc}
\usepackage[english]{babel}
\usepackage{xcolor}

\usepackage{amsmath, amssymb, bm}

\usepackage{graphicx}
\usepackage{booktabs}
\usepackage{multirow}
\usepackage{makecell}
\usepackage{enumitem}

\usepackage[margin=1.5in]{geometry}
\makeatletter
\renewcommand\footnoterule{%
  \vspace{0.7em}   
  \hrule width 0.4\columnwidth height 0.4pt
  \kern 0.3em      
}
\makeatother
\usepackage[hang]{footmisc}
\usepackage{float}
\usepackage{tabularx}
\usepackage{natbib}
\usepackage{hyperref}
\hypersetup{
    colorlinks=false,    
    pdfborderstyle={},   
    pdfborder={0 0 0},  
    linkcolor=black,
    citecolor=black,
    filecolor=black,
    urlcolor=black
}

\begin{document}

\definecolor{frontMinMax}{HTML}{2E8B57}
\definecolor{frontMinMin}{HTML}{8C564B}
\definecolor{frontMaxMax}{HTML}{E07B39}
\definecolor{frontMaxMin}{HTML}{7B2CBF}

\title{\textbf{Preference-Based Constraint Handling for Multi-Objective Design-Decision Optimisation: Acceptability Relaxation and Feasibility Guidance}}

\renewcommand{\thefootnote}{\fnsymbol{footnote}}

\author{
E. \"Ozer$^{1}$, A.R.M. Wolfert$^{2}$\thanks{Corresponding author: A.R.M. (Rogier) Wolfert, \texttt{a.r.m.wolfert@tudelft.nl}}\\[0.5em]
\small$^{1}$ Decision Science \& Engineering, Boskalis, The Netherlands\\[-0.3em]
\small$^{2}$ Department of Algorithmics, Faculty of Mathematics and Computer Science\\[-0.3em]
\small Delft University of Technology, The Netherlands
}
\date{}
\maketitle
\renewcommand{\thefootnote}{\arabic{footnote}}
\vspace{-2.0em}

\begin{abstract}

\noindent
Multi-objective design-decision optimisation (MODO) problems often contain both negotiable acceptability limits and non-negotiable feasibility constraints. Hard constraints define the feasible design space and cannot be violated, whereas negotiable limits represent initially acceptable performance levels that actors may reconsider through deliberation and compromise. Conventional optimisation methods commonly handle constraint relaxation through deviations, penalties, utility values, desirability scores, or violation measures. However, relaxing an acceptability limit is fundamentally a preference judgement about how actors value departures from the original requirement. Trade-offs expressed solely through technical measures cannot therefore be assumed to represent meaningful preference trade-offs and may yield recommendations that are neither preference-consistent nor decision-valid.

\noindent
This paper introduces a preference-based constraint-handling framework for MODO that explicitly distinguishes negotiable acceptability limits, which have preference and decision meaning, from non-negotiable feasibility constraints, which remain fixed and permit only computational search guidance. For negotiable limits, preference functions represent how actor preference changes as performance departs from the initially acceptable level, allowing both the preferred design and the acceptable degree of relaxation to be determined endogenously through the same preference-aggregation process. Two mathematically distinct formulations are developed: Approach~1 represents each relaxation explicitly as an additional objective and search variable, whereas Approach~2 incorporates relaxation directly into the preference function of the affected criterion. When constructed consistently, both approaches yield the same solution that best fits the common purpose, while differing in modelling complexity, interpretability, and computational effort.

\noindent
The framework addresses two complementary acceptability-relaxation situations. In over-constrained stalemate problems, where no design satisfies all negotiable limits, it identifies the preference-consistent relaxations required to recover an acceptable solution. In non-stalemate problems, where feasible designs already exist, it determines whether modest relaxation can reveal a substantially better design and whether the resulting gain in aggregate actor preference justifies the departure from the original limits.

\noindent
For hard feasibility constraints, the same mathematical construction is assigned a fundamentally different role. Hard constraints remain fixed, and infeasible designs remain inadmissible. A feasibility-guidance function instead provides the IMAP-IGS search with graded computational information about scaled constraint violations, enabling it to distinguish infeasible candidates and progress towards small, fragmented, or otherwise difficult-to-reach feasible regions. Every feasible design takes precedence over every infeasible design, and the final recommendation is selected exclusively from the feasible set using the A-FINE aggregate preference score.

\noindent
The resulting framework provides an integrated preference-based approach to constraint handling in MODO: it enables recovery from over-constrained stalemate, identifies justified improvement opportunities beyond negotiable limits in non-stalemate problems, and guides the search towards hard feasibility boundaries in highly constrained design spaces. It thereby integrates acceptability negotiation, feasibility guidance, and design optimisation without conflating their distinct meanings, while preserving preference consistency and decision validity.

\end{abstract}

\section{Introduction}

Constraint handling is a fundamental challenge in engineering design and managerial decision-making whenever performance limits cannot simultaneously be satisfied or may be reconsidered to achieve more preferred outcomes. This section reviews existing constraint-handling approaches, motivates a preference-based framework for negotiable constraint relaxation and non-negotiable constraint guidance, and introduces the underlying problem statement.

\subsection{Current Approaches to Constraint Handling and Their Limitations}
\label{sec:limitations}

The treatment of performance limits and constraint handling has long been studied in engineering design optimisation, systems engineering, operations research, and decision analysis \citep{Ehrgott2005,Blanchard2021,Hillier2021,Martins2022}. Numerous approaches have been proposed for situations in which one or more performance requirements cannot be satisfied simultaneously or the feasible region is difficult to reach.

Goal and elastic programming represent performance limits as goals and introduce deviation variables that are subsequently minimised according to predetermined weights, priorities, or minimax criteria \citep{Charnes1961,Hillier2021}. Fuzzy and desirability-based approaches replace rigid limits with membership functions and seek solutions that maximise aggregate satisfaction levels \citep{Bellman1970,Derringer1980}. Multi-attribute utility methods assign utility functions to individual criteria and combine them through additive, multiplicative, or minimum-based aggregation rules \citep{Keeney1976}. Reference-point and bargaining approaches evaluate alternatives relative to aspiration levels or disagreement points \citep{Nash1950}, whereas constrained multi-objective optimisation formulations treat constraint violations as additional objectives and identify Pareto-efficient trade-off solutions \citep{Ehrgott2005,Martins2022}.

These approaches are well established in engineering optimisation and decision support, but they primarily provide mathematical mechanisms for managing infeasibility, deviations, and compromise. Their ability to produce a computationally or mathematically valid solution does not establish that the underlying trade-offs reflect actor preferences or that the resulting recommendation is preference-consistent and decision-valid. They therefore remain limited when constraint handling is treated as a multi-objective design-decision optimisation (MODO) problem. In other words, mathematical tractability and optimisation validity do not, by themselves, establish decision validity or resolve the underlying design-decision problem. This limitation may also impede practical adoption in engineering-design and project-management organisations, where decision makers require transparent and interpretable justification for why particular constraints are relaxed, which trade-offs are accepted, and how the resulting recommendation reflects the preferences of the actors involved. Without such a design-decision orientation, technically sophisticated constraint-handling methods may leave the underlying acceptability question unresolved and provide limited support for reaching a clear, collectively justified decision. More specifically, five interrelated limitations prevent these approaches from fully addressing constraint handling as a human-centred, preference-based MODO problem:

\begin{enumerate}
    \item \textbf{Failure to distinguish acceptability from feasibility.} Conventional methods often treat negotiable acceptability limits and non-negotiable feasibility constraints as comparable mathematical restrictions, despite their fundamentally different meanings. Negotiable limits represent socially, institutionally, or contractually established judgements and may be reconsidered through actor deliberation and preference trade-offs. Hard constraints represent physical, technical, legal, or safety boundaries that remain fixed and permit only computational guidance towards feasibility. Although ODESYS in FIVES explicitly recognises this distinction \citep{Wolfert2026b}, conventional constraint-handling methods do not generally preserve it.

    \item \textbf{Representation-dependent and potentially invalid trade-offs.} Recommendations may depend on whether departures from specified limits are represented as absolute or relative deviations, penalties, utilities, desirability scores, membership grades, or violation measures. Because these representations possess different measurement properties, their numerical aggregation does not automatically constitute a meaningful preference trade-off \citep{Barzilai2022, Wolfert2026a}.

    \item \textbf{Hidden modelling assumptions.} Penalty weights, aspiration levels, priority structures, utility-function shapes, membership functions, tolerance bands, and aggregation operators can substantially influence the result, yet are often selected or calibrated by analysts rather than elicited directly from actor preferences \citep{Wolfert2026b, Timp2026}.

    \item \textbf{Lack of a unique recommendation.} Pareto-based and bargaining-oriented methods \citep{Nash1950} may identify admissible compromises without determining which compromise is preferred. The outcome is therefore often a set of alternatives, leaving the final choice to an additional decision process outside the optimisation \citep{Martins2022,Wolfert2026b,OzerWolfert2026,Timp2026}.

    \item \textbf{Insufficient guidance towards feasibility.} In highly constrained problems, the feasible region may be small, fragmented, or difficult to reach. Binary rejection or equal penalisation of infeasible candidates provides no information about their relative progress towards feasibility, causing population-based solvers to stagnate or fail to identify a feasible design. Graded constraint-violation information can guide the search, but must remain computational and separate from the final preference-based ranking of feasible designs \citep{Wolfert2026b,Timp2026}.
\end{enumerate}

More fundamentally, existing constraint-handling approaches often fail to distinguish between the decision meaning of negotiable acceptability limits and the computational role of hard feasibility constraints. Relaxing a negotiable acceptability limit is a preference construct: the relevant question is how actors value departures from an originally acceptable performance level. By contrast, violation of a hard feasibility constraint has no decision meaning and can never be justified by preference. Information about such violations may nevertheless serve as a computational search signal for guiding an optimisation algorithm towards the feasible region.

For negotiable acceptability limits, relaxation is commonly represented through deviations, penalties, utilities, desirability measures, distances, or violation scores that must somehow be balanced against one another. Yet the actual decision concerns the willingness of actors to accept departures from originally acceptable performance levels. The quantity of decision interest is therefore not the physical magnitude of the relaxation itself, but the preference associated with that relaxation.

From a measurement-theoretical perspective, a design-decision recommendation is valid only if trade-offs occur in a commensurable preference space rather than directly on physical quantities, deviations, utility values, or violation measures. Recent work on preference-based optimisation and aggregation has argued that preference consistency, commensurability, and preference preservation \citep{Barzilai2005,Barzilai2010} are necessary conditions for meaningful design-decision support \citep{Wolfert2023,Wolfert2026a,Wolfert2026b,Timp2026}. Existing constraint-relaxation approaches generally provide no such guarantee. Consequently, the resulting recommendation may be optimisation-valid but not necessarily decision-valid.

From this perspective, preference-based constraint handling requires a strict distinction between acceptability relaxation and feasibility guidance. Negotiable performance limits express actor preferences regarding acceptable levels of system performance and should therefore be modelled through preference functions that quantify the acceptability of departures from those limits. Both the preferred design and the acceptable degree of relaxation can then emerge from a single preference-consistent design-decision process. Hard feasibility constraints, by contrast, remain fixed and cannot be relaxed through preference trade-offs. In highly constrained optimisation problems, scaled constraint violations may nevertheless provide computational guidance that helps the solver approach the feasible region, provided that this guidance remains separate from the preference-based ranking and final recommendation.

Accordingly, this paper develops an integrated preference-based constraint-handling framework for MODO with three complementary applications: recovering acceptable solutions in over-constrained stalemate problems, identifying justified improvement opportunities beyond negotiable limits in non-stalemate problems, and guiding the IMAP-IGS search towards hard feasibility boundaries in highly constrained design spaces.

The central distinction is therefore that acceptability constraints require preference-consistent decision modelling, whereas feasibility constraints require computational search guidance without decision meaning. By preserving this distinction, the framework integrates acceptability negotiation, design optimisation, and feasibility guidance without allowing hard-constraint violations to influence the final preference-based recommendation.
\subsection{Problem Formulation}

The limitations identified in the preceding subsection indicate that constraint handling requires a strict distinction between negotiable acceptability limits and hard feasibility constraints. Relaxation of a negotiable limit should not be treated merely as a technical violation-minimisation problem, but as a preference-based design-decision problem. By contrast, violations of hard constraints cannot be accepted through preference trade-offs, although information about such violations may be used to guide the optimisation search. These principles are illustrated using a grillage-design case study for a marine contractor.

Each candidate design is characterised by its grillage type and steel quantity and is evaluated in terms of investment cost \(c\), construction time \(t\), and structural safety margin \(g\). The safety requirement \(g\geq g_{\min}\) is a hard feasibility constraint that must be satisfied by every admissible design. Cost and construction time, by contrast, are subject to negotiable acceptability limits represented by the budget limit \(B\) and project-duration limit \(T\), respectively.

For the negotiable limits, the central question is not simply how far a design exceeds a limit, but how acceptable that departure is to the actors involved. Exceeding the budget or schedule limit represents a loss of actor preference that must be considered alongside the preference gains achieved elsewhere in the design. Acceptability-constraint relaxation therefore constitutes a design-decision problem rather than a purely technical feasibility adjustment.

Two acceptability-relaxation situations are considered. The first is an \emph{over-constrained stalemate case}, in which no design simultaneously satisfies all negotiable limits. At least one acceptability limit must then be relaxed to recover an acceptable solution, and the decision problem becomes determining which limits may be relaxed and by how much. Two alternative formulations are developed. Approach~1 introduces an explicit relaxation objective for each negotiable constraint, whereas Approach~2 incorporates relaxation directly into the criterion-specific preference functions. When the underlying preference functions and weights are constructed consistently, both approaches produce the same best-fit solution.

The second situation is a \emph{non-stalemate case}, in which one or more designs already satisfy all negotiable limits. Relaxation is then not required to obtain an acceptable design. Nevertheless, a modest departure from the original budget or schedule limit may reveal designs with substantially higher aggregate actor preference than the best design within the original acceptability limits. The framework therefore determines whether the resulting preference gain justifies the relaxation and, if so, the acceptable extent of that relaxation.

A third application concerns feasibility-constraint handling in highly constrained design spaces. The hard safety requirement remains fixed and cannot be relaxed or traded against actor preferences. However, if the feasible region is small, fragmented, or difficult to reach, a binary feasible-infeasible classification may provide insufficient information for the search algorithm. A feasibility-guidance function is therefore used to distinguish infeasible candidates according to their scaled constraint violations and guide the search towards the feasible region. These guidance values are computational search signals without preference or decision meaning. Every feasible design takes precedence over every infeasible design, and the final recommendation is always selected from the feasible set.

The problem addressed in this paper is therefore threefold: to recover an acceptable solution when negotiable limits create a stalemate, to identify justified improvement opportunities beyond negotiable limits when acceptable designs already exist, and to guide the IMAP-IGS search towards hard feasibility boundaries in highly constrained design spaces. In the first two applications, the preferred design and acceptable degree of relaxation emerge jointly from the same preference-consistent aggregation process. In the third, feasibility guidance supports the search without relaxing hard constraints or influencing the final preference-based ranking of feasible designs.

\section{Constraint Relaxation Through Preference Optimisation}
\label{sec:method}

Constraint relaxation becomes necessary whenever negotiable performance limits prevent the identification of a satisfactory design. Rather than treating relaxation as a technical violation-minimisation problem, the proposed approach formulates it as a preference-based design--decision problem in which both the design and the acceptable degree of relaxation are determined simultaneously through preference optimisation. Two equivalent formulations are presented, differing only in the way relaxation preferences are modelled and incorporated into the optimisation process.

\subsection{Core idea, notation, and formulation choice}

Let \(\mathbf{x}\) denote a design characterised by its cost \(c(\mathbf{x})\), construction time \(t(\mathbf{x})\), and safety performance \(g(\mathbf{x})\). The cost and time limits, denoted by \(B\) and \(T\), respectively, are negotiable constraints. For a given design, the minimum relaxations required to satisfy these limits are

\begin{equation}
\label{eq:required-relaxations}
r_c(\mathbf{x})
=
\max\left\{0,c(\mathbf{x})-B\right\},
\qquad
r_t(\mathbf{x})
=
\max\left\{0,t(\mathbf{x})-T\right\}.
\end{equation}

Thus, \(r_c(\mathbf{x})=0\) and \(r_t(\mathbf{x})=0\) when the original cost and time limits are satisfied, whereas positive values quantify the minimum relaxations required to render the design acceptable with respect to these limits. By contrast, the safety requirement \(g(\mathbf{x})\geq g_{\min}\) is a hard constraint and must be satisfied by every admissible design.

Throughout, \(k\) indexes actors and \(i\) indexes performance criteria. The relaxation formulations below are stated at the criterion level; actor-specific preference terms are not expanded explicitly. Because acceptability relaxation requires a judgement about the extent to which departures from the original limits remain acceptable, it is treated as a design-decision problem rather than merely as a technical feasibility adjustment. For each criterion \(i\), a preference function \(p_i\) \citep{Barzilai2005,Barzilai2010} maps criterion performance to a preference score on a scale from \(0\) to \(100\) \citep{Wolfert2026b}, thereby representing how preference changes as performance departs from the original acceptable limit. The resulting preference scores are aggregated through the affine IMAP aggregation procedure \citep{Wolfert2023,Wolfert2026a,Wolfert2026b}, and the IGS-GA solver maximises the aggregate preference. The safety threshold and the maximum permissible relaxation of each negotiable limit remain hard bounds. Consequently, the preferred design and the acceptable degrees of relaxation are determined simultaneously through the same preference-based optimisation process \citep{Wolfert2026b}.

This preference model can be formulated through two mathematically distinct but equivalent approaches. In \emph{Approach~1, the explicit-relaxation formulation}, each relaxation is represented as an additional objective and search variable. In \emph{Approach~2, the integrated-relaxation formulation}, relaxation is incorporated directly into the preference function of the corresponding performance criterion. When the preference functions and weights are constructed consistently, both approaches produce the same best-fit design, as demonstrated in Section~\ref{sec:equiv}. Their difference therefore concerns model construction and computational bookkeeping rather than the underlying feasible set or decision logic.

For a problem with \(n\) negotiable criteria, the practical differences between the two approaches are summarised in Table~\ref{tab:select}.

\begin{table}[htbp]
\centering

\label{tab:select}
\small
\setlength{\tabcolsep}{5pt}
\renewcommand{\arraystretch}{1.15}

\begin{tabularx}{\textwidth}{
    >{\raggedright\arraybackslash}p{0.25\textwidth}
    >{\raggedright\arraybackslash}X
    >{\raggedright\arraybackslash}X
}
\toprule
\textbf{Model aspect}
&
\textbf{Approach~1: explicit relaxation}
&
\textbf{Approach~2: integrated relaxation}
\\
\midrule
Preference functions
&
\(2n\)
&
\(n\)
\\
Weights to specify
&
\(2n\)
&
\(n\)
\\
Search variables
&
Design variables and \(n\) relaxation variables
&
Design variables only
\\
Representation of relaxation
&
Explicit relaxation objective
&
Implied by the selected design
\\
\bottomrule
\end{tabularx}
\caption{Comparison of the two relaxation approaches for \(n\) negotiable criteria.}
\end{table}

As a practical rule, Approach~1 is generally preferable for problems with only two or three negotiable criteria. Each relaxation then has its own linear preference function and weight, allowing the decision maker to specify directly how an overrun of one limit should be valued relative to an overrun of another. Because the relaxation amounts are explicit objectives, they can also be reported, interpreted, and discussed separately.

Approach~2 is generally more efficient for problems with four or more negotiable criteria. Applying Approach~1 would require twice as many preference functions and weights, while each explicit relaxation would introduce an additional search variable. Approach~2 retains one preference function per criterion and searches only over the original design variables. The preference trade-off associated with relaxation is then specified within each criterion-level preference function through its anchors and the preference score assigned to the original limit.

The choice between the two approaches therefore affects model complexity, interpretability, and computational effort, but not the applicable hard bounds. In both approaches, the safety threshold and the maximum permissible relaxations remain hard constraints, while the acceptable degree of relaxation emerges endogenously from the preference-based optimisation.

\subsection{Approach~1: explicit relaxation objectives}

Each criterion retains its own preference function, and each relaxation becomes a separate performance metric with its own preference function. With two criteria, the decision variables are the design $\mathbf{x}$ and the relaxation variables $\Delta_c$ and $\Delta_t$:
\begin{equation}
\max_{\mathbf{x},\,\Delta_c,\,\Delta_t}\;
\sum_{i\in\{c,t\}}\Bigl[\,w_i\,p_i\bigl(y_i(\mathbf{x})\bigr)+w_{\Delta i}\,p_{\Delta i}(\Delta_i)\Bigr]
\qquad\text{s.t.}\qquad
\begin{aligned}
&y_c(\mathbf{x})\le B+\Delta_c\le B+\bar\Delta_c,\\
&y_t(\mathbf{x})\le T+\Delta_t\le T+\bar\Delta_t,\\
&g(\mathbf{x})\ge g_{\min},\qquad \Delta_i\ge 0,
\end{aligned}
\label{eq:part1}
\end{equation}
where $y_c=c$, $y_t=t$, the weights satisfy $\sum_i(w_i+w_{\Delta i})=1$, and $\bar\Delta_i$ is the maximum relaxation of each limit.

Each criterion's preference function is linear over the performance range of the design space, scoring 100 at the best value and 0 at the worst. Each relaxation preference function is linear, scoring 100 at zero relaxation and 0 at the maximum relaxation:
\begin{equation*}
p_{\Delta i}(\Delta)=100\Bigl(1-\frac{\Delta}{\bar\Delta_i}\Bigr).
\end{equation*}
The study uses weights $(w_c,w_t,w_{\Delta c},w_{\Delta t})=(0.2,0.2,0.3,0.3)$. Because $p_{\Delta i}$ decreases, the best relaxation for a fixed design is the smallest admissible one, $\Delta_i^\star=r_i(\mathbf{x})$. The solver determines this relaxation jointly with the design.

\begin{figure}[htbp]
\centering
\includegraphics[width=0.9\textwidth]{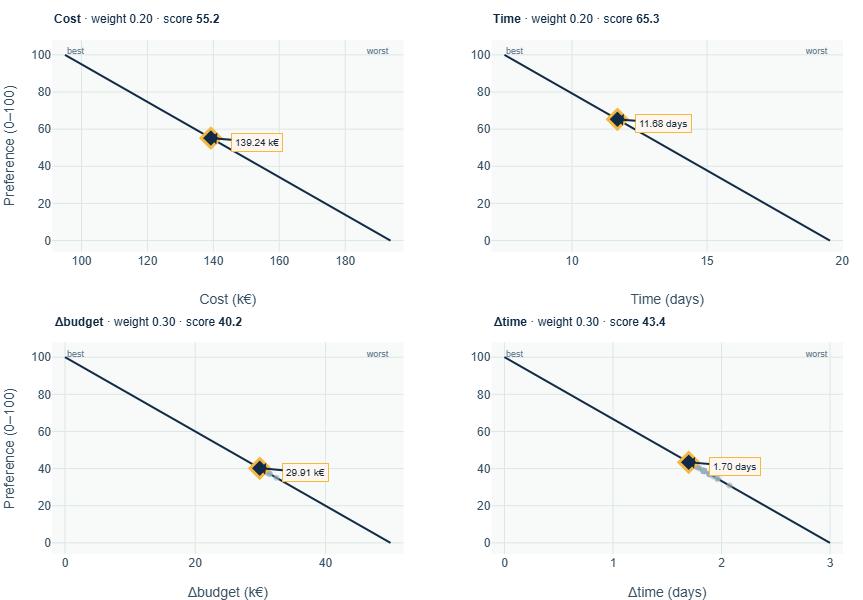}
\caption{Preference curves for the selected Design C. The criterion curves score design performance, while the relaxation curves score departures from the budget and schedule limits.}
\label{fig:part1-curves}
\end{figure}

\begin{figure}[htbp]
\centering
\includegraphics[width=0.9\textwidth]{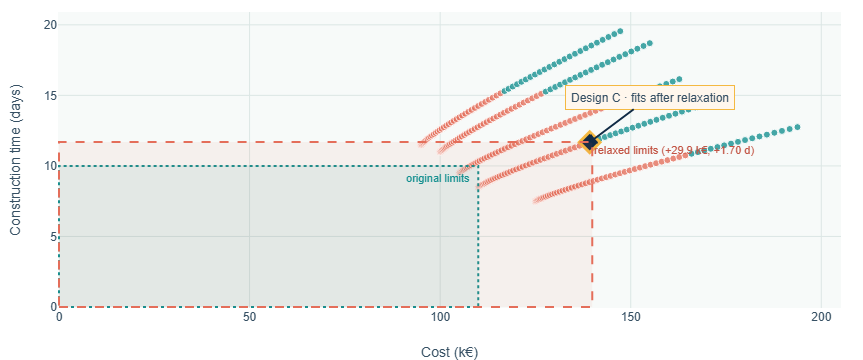}
\caption{Original and relaxed acceptable regions for Approach~1. Candidate designs are shown relative to both regions, with the selected Design C marked by a diamond.}
\label{fig:part1-region}
\end{figure}

The original budget in this example case is 110\,k\texteuro{} and the time limit is 10 days, but no design satisfies both limits. The Approach~1 optimisation selects Design C, with a cost of 139.24\,k\texteuro{} and a construction time of 11.68 days. This result requires raising both the budget and the time limit. Figure~\ref{fig:part1-curves} shows the preference curves for Design C, and Figure~\ref{fig:part1-region} shows its position relative to the original and relaxed acceptable regions.
For further numerical analysis, readers are referred to the Data Availability Statement.

\subsection{Approach~2: modified preference curves}
\label{sec:part2}

Approach~2 uses one curve per criterion and incorporates the relaxation into that curve \citep{Wolfert2026b}. Each curve has three anchors: the best value $a_i$ (score 100), the acceptable limit $L_i$ (limit score $s_i$), and the largest relaxed value $U_i=L_i+\bar\Delta_i$ (score 0). The curve is linear from the best value to the limit, then falls steeply to 0 at the largest relaxed value:
\begin{equation*}
q_i(y)=
\begin{cases}
100, & y\le a_i,\\[3pt]
100-(100-s_i)\,\dfrac{y-a_i}{L_i-a_i}, & a_i<y\le L_i,\\[9pt]
s_i\,\dfrac{e^{-\kappa u}-e^{-\kappa}}{1-e^{-\kappa}}, & L_i<y<U_i,\quad u=\dfrac{y-L_i}{U_i-L_i},\\[9pt]
0, & y\ge U_i.
\end{cases}
\end{equation*}
The limit score $s_i$ represents the decision maker's preference at the acceptable limit. The steepness $\kappa$ determines how quickly the score falls beyond that limit. The normalisation in the exponential term ensures that $q_i(U_i)=0$ exactly.

\begin{figure}[htbp]
\centering
\includegraphics[width=0.9\textwidth]{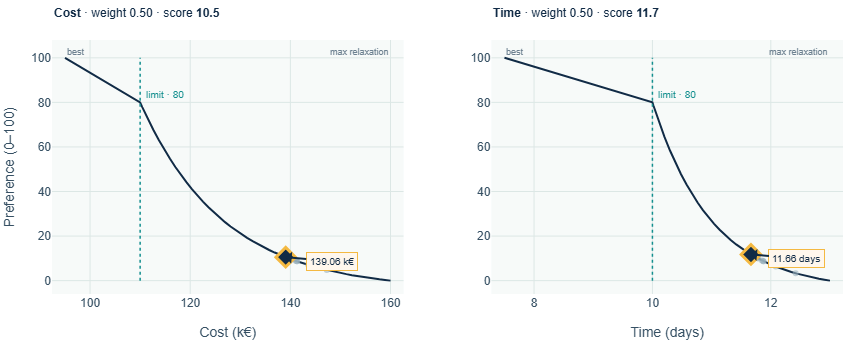}
\caption{Preference curves used in the two-objective study. Preference declines beyond each negotiable limit, and the selected design lies beyond both original limits.}
\label{fig:curves}
\end{figure}

Approach~2 then solves
\begin{equation}
\max_{\mathbf{x}}\; v_c\,q_c\bigl(c(\mathbf{x})\bigr)+v_t\,q_t\bigl(t(\mathbf{x})\bigr)
\qquad\text{s.t.}\qquad
c(\mathbf{x})\le U_c,\quad t(\mathbf{x})\le U_t,\quad g(\mathbf{x})\ge g_{\min},
\label{eq:part2}
\end{equation}
with weights $v_c+v_t=1$. The relaxation is not a decision variable. It is read off the selected design as $r_i(\mathbf{x}^\star)=\max\{0,\,y_i(\mathbf{x}^\star)-L_i\}$. The case-study implementation uses $s_i=80$, $\kappa=3$ and $v_c=v_t=0.5$, which produce the curves in Figure~\ref{fig:curves}. In this stalemate case, no design satisfies the original budget of 110\,k\texteuro{} and time limit of 10 days simultaneously. The optimisation selects a design costing 139.24\,k\texteuro{} and taking 11.68 days as the best under this preference model, requiring the budget to increase by 29.24\,k\texteuro{} and the time limit by 1.68 days.
For further numerical analysis, readers are referred to the Data Availability Statement.

\subsection{Why both approaches give the same design-decision outcome}
\label{sec:equiv}

Let $v_i=w_i+w_{\Delta i}$, so that $v_c+v_t=1$, and define the combined curve
\begin{equation*}
q_i(y)=\frac{w_i\,p_i(y)+w_{\Delta i}\,p_{\Delta i}\bigl(\max\{0,\,y-L_i\}\bigr)}{v_i}.
\end{equation*}
For every design $\mathbf{x}$ within the upper bounds $U_i$, the Approach~1 objective at its optimal relaxation equals $\sum_i v_i\,q_i\bigl(y_i(\mathbf{x})\bigr)$. The two formulations therefore have the same optimal design.

\textbf{Proof.} Fix $\mathbf{x}$. Approach~1 constraints require $\Delta_i\ge r_i(\mathbf{x})$, and $p_{\Delta i}$ decreases, so the optimal relaxation is $\Delta_i^\star=\max\{0,\,y_i(\mathbf{x})-L_i\}$. Substituting it into \eqref{eq:part1} gives
\begin{equation*}
\sum_i\Bigl[w_i\,p_i(y_i)+w_{\Delta i}\,p_{\Delta i}\bigl(\max\{0,\,y_i-L_i\}\bigr)\Bigr]=\sum_i v_i\,q_i(y_i),
\end{equation*}
which is the objective of \eqref{eq:part2}. The feasible designs are the same in both models, so the maximisers coincide.

The equivalence is exact for the combined curves $q_i$. The curves in Approach~2 (Section~\ref{sec:part2}) are specified directly and encode their own trade-off. If they are not derived from the Approach~1 curves, the two models can return different designs. To ensure consistency, derive the Approach~2 curves from the Approach~1 curves as above and evaluate both formulations on the same case. With the default four-objective weights, the implied two-objective weights are $v_c=v_t=0.5$, matching the default weights in the two-objective study.

\section{Relaxation Can Reveal an Even Better Optimum}
\label{sec:opportunity}

The previous sections addressed over-constrained stalemate problems in which no design satisfies all specified limits. A second situation arises when feasible designs already exist, but relaxing one or more negotiable limits may reveal further opportunities for improvement. The limits are negotiable, but the decision maker treats all performance beyond them as unacceptable, so the region beyond them is not searched. Relaxing a limit slightly can open that region, which may contain a design that is clearly better than the best feasible one. The relaxation must still be judged against the feasible alternatives, because a larger region does not guarantee a better design.

We illustrate this case with a similar example. The limits are a budget of 110\,k\texteuro{} and a time limit of 10 days; the safety threshold remains a hard constraint. The preference curves are those of Approach~2 in Section~\ref{sec:part2}, with limit score $s_i=40$ and steepness $\kappa=1$. Each curve falls linearly from 100 at its best value to 40 at its limit, then decays to 0 at the maximum relaxation of 30\,k\texteuro{} and 3 days. The weights are $0.6$ for cost and $0.4$ for time. The study proceeds in three stages:
\begin{enumerate}
\item Find the best design inside the original limits, the Original IMAP Optimum.
\item Rerun with the limits widened to the maximum relaxation. The original optimum is seeded into the initial population, so a relaxed design must outrank it to be selected.
\item Score both optima against a common reference pool containing all designs that satisfy the hard safety constraint and the maximum relaxed limits, so that their scores are comparable.
\end{enumerate}

The two figures below, Figure~\ref{fig:opp-region} and Figure~\ref{fig:opp-map}, demonstrate the original and relaxed acceptable regions and where a better optimum can lie.

\begin{figure}[htbp]
\centering
\includegraphics[width=0.9\textwidth]{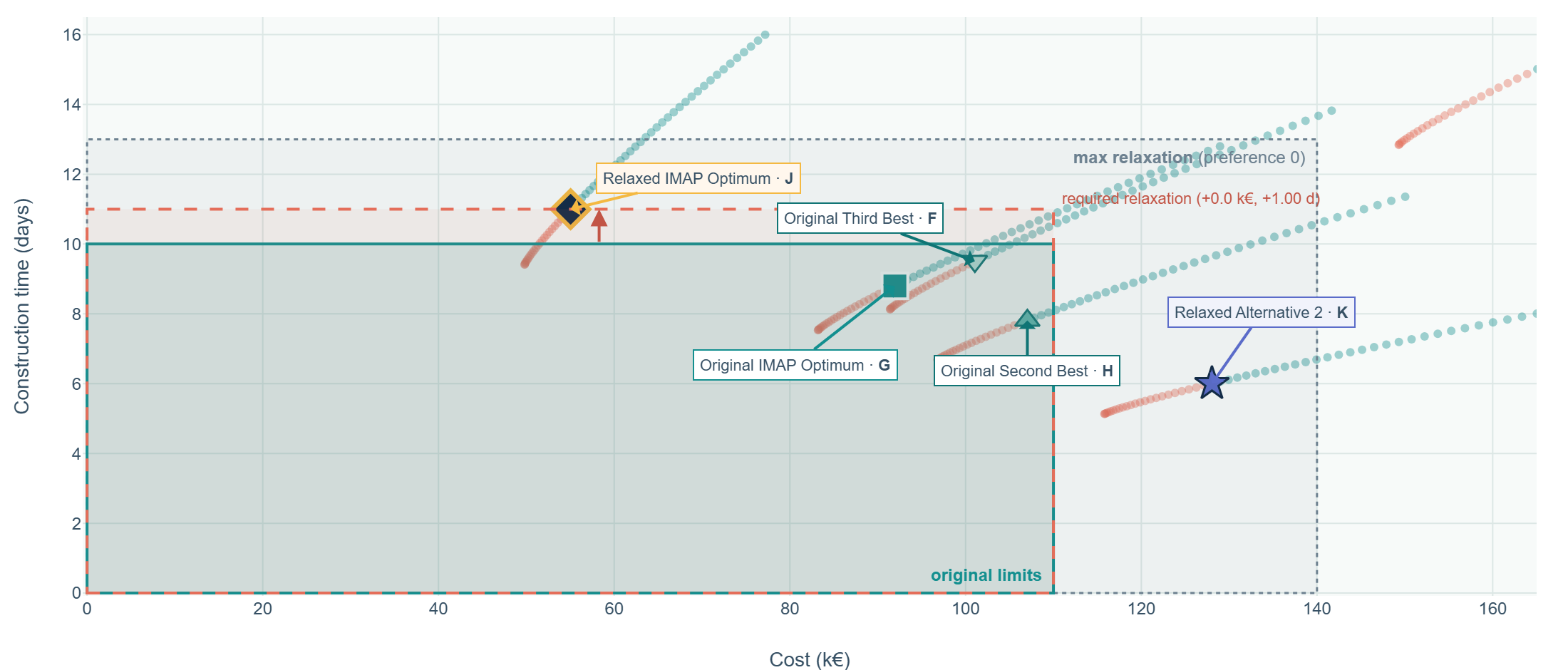}
\caption{Original, required-relaxation, and maximum-relaxation regions for the non-stalemate study. Sampled designs and labelled alternatives show how expanding the negotiable limits changes the admissible region.}
\label{fig:opp-region}
\end{figure}

\begin{figure}[htbp]
\centering
\includegraphics[width=0.9\textwidth]{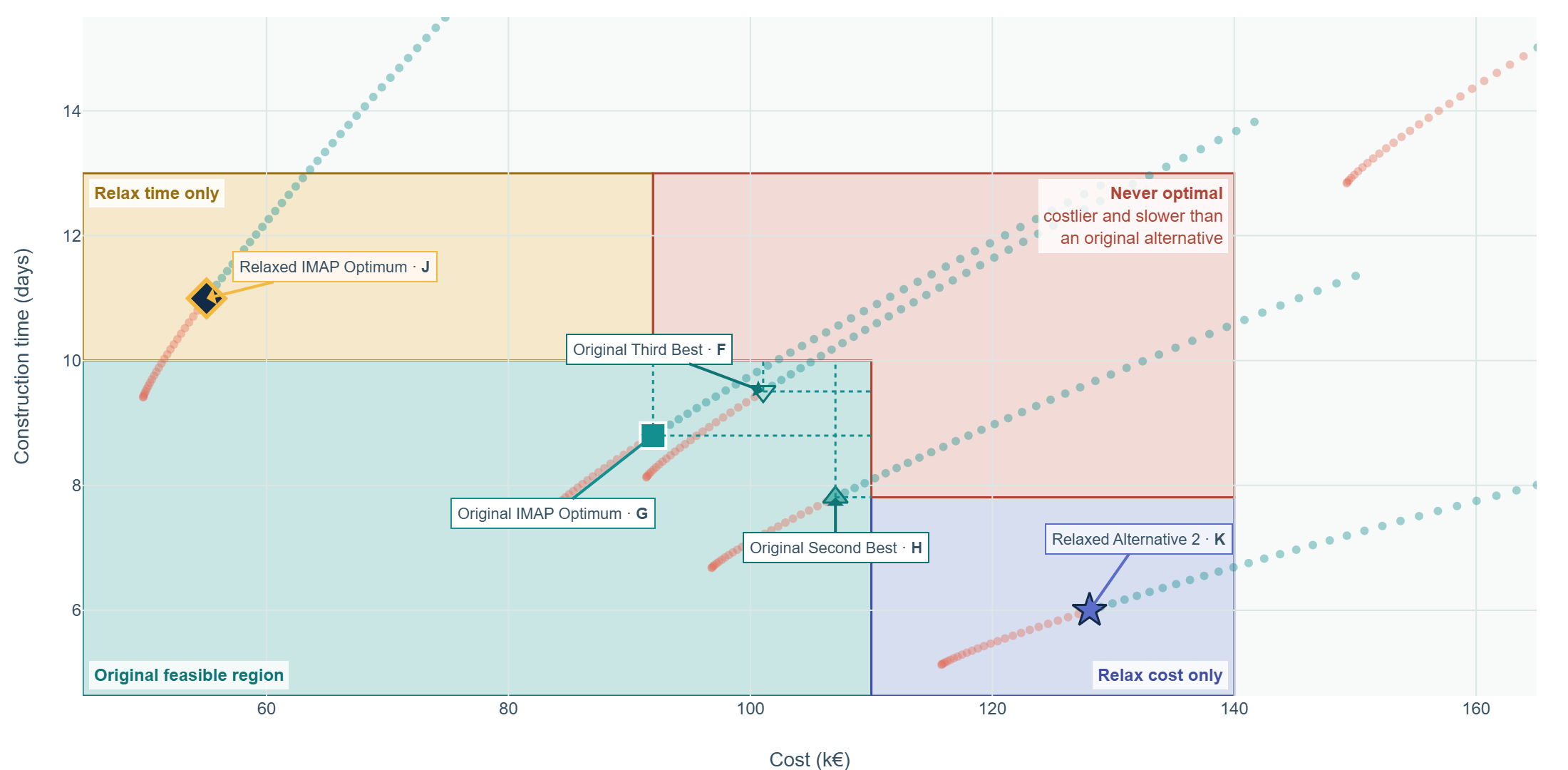}
\caption{Cost-time regions that distinguish designs dominated by an original feasible alternative from candidates that may become preferred after relaxation.}
\label{fig:opp-map}
\end{figure}

Figure~\ref{fig:opp-map} divides the cost-time plane into four regions. A design in the red region is at least as costly and as slow as an original alternative, which therefore dominates it. Such a design cannot score higher than that alternative and cannot be optimal. Only the amber and blue regions can contain designs that are not dominated by an original alternative. In the study, J lies in the amber region and K in the blue region. The size of the gain depends on the preference settings and weights. The scores of 100 and 0 describe this set only; they are not absolute measures of quality.

Table~\ref{tab:alternatives} lists the five named alternatives. The IMAP score is rescaled across these alternatives, so the values are relative to this set.

\begin{table}[htbp]
\centering
\small
\resizebox{\textwidth}{!}{%
\begin{tabular}{llccccl}
\toprule
Alternative & Design & Steel (kg/m$^2$) & Cost pref. & Time pref. & IMAP score & Required relaxation \\
\midrule
Original IMAP Optimum & G Composite frame & 43.94 & 58.0 & 54.8 & 64.1 & No \\
Original Second Best & H Modular frame & 44.00 & 43.0 & 67.0 & 46.3 & No \\
Original Third Best & F Standard frame & 44.00 & 48.9 & 46.1 & 9.4 & No \\
\midrule
Relaxed IMAP Optimum & J Precast economy & 43.94 & 94.8 & 22.1 & 100.0 & Yes, time +1.00 day \\
Relaxed Alternative 2 & K Prefab rapid & 44.00 & 11.4 & 89.3 & 0.0 & Yes, budget +18.0\,k\texteuro{} \\
\bottomrule
\end{tabular}%
}
\caption{Original versus relaxed alternatives: preference scores and limit relaxations.}
\label{tab:alternatives}
\end{table}

The best design in the original space is G, at 92\,k\texteuro{} and 8.8 days, which has an IMAP score of 64.1 when compared to newly discovered alternatives. Relaxing the time limit by one day admits J at 55\,k\texteuro{} and 11.0 days. J is one day over the time limit, but it costs 37\,k\texteuro{} less than G and receives the highest score in the set. The study therefore recommends the time relaxation, provided the decision maker accepts the extra day. Relaxing the budget by 18\,k\texteuro{} admits K at 128\,k\texteuro{}. K is faster than the feasible alternatives shown, but it receives the lowest score in the set. A relaxed design is therefore not automatically better, and each one must be compared with the feasible alternatives.
For further numerical analysis, readers are referred to the Data Availability Statement.

\section{Feasibility Guidance for IMAP--IGS Metaheuristics}
\label{sec:guidance}

Sections~\ref{sec:method} and~\ref{sec:opportunity} addressed preference-based acceptability-constraint handling in two complementary situations: over-constrained stalemate problems, in which no design satisfies all negotiable limits, and non-stalemate problems, in which feasible designs exist but limited relaxation may reveal better solutions. As an additional application, an analogous mathematical construction can support feasibility-constraint handling for hard constraints, such as the safety threshold \(g(\mathbf{x})\geq g_{\min}\).

In this application, the purpose is fundamentally different. A hard constraint is never relaxed, and infeasible designs remain inadmissible. Rather than representing actor preference, a feasibility-guidance function provides the IGS solver with graded computational information about the progress of an infeasible design towards the feasible region. This function may be constructed analogously to the preference function used in Approach~2 \citep{Wolfert2026b}, but its values have no preference or decision meaning.

Such guidance is particularly important in highly constrained problems, where the feasible region may be small, fragmented, or difficult to reach, and a binary feasible-infeasible distinction may provide insufficient information for effective search \citep{Wolfert2026b}. The feasibility-guidance function differentiates infeasible designs according to their progress towards feasibility while preserving the hard feasibility boundary. It thereby supports solver convergence and increases the likelihood of finding a feasible solution without allowing the hard constraint to influence the preference-based ranking of feasible designs.

Table~\ref{tab:guidance} compares preference-based handling of negotiable acceptability limits with computational guidance for hard feasibility constraints. Although both applications may use analogous mathematical constructions, they answer fundamentally different questions.

\begin{table}[htbp]
\centering
\small
\label{tab:guidance}

\begin{tabularx}{\textwidth}{>{\raggedright\arraybackslash}p{3.4cm}XX}
\toprule
 & Negotiable acceptability limits (Sections~\ref{sec:method}--\ref{sec:opportunity}) & Hard feasibility constraints (this section) \\
\midrule
Purpose & Determine whether, and to what extent, a limit may be relaxed & Guide the solver towards the feasible region \\
Constraint after the search & May be relaxed when the preference gain justifies it & Remains fixed and must be satisfied \\
Meaning beyond the limit & Actor preference for departures from the original acceptable limit & Computational signal indicating progress towards feasibility \\
Recommendation beyond the limit & Permissible within the maximum relaxation when justified by aggregate preference & Never permissible \\
Role in the final decision & Contributes to the preference-based ranking and recommendation & Does not enter the preference-based ranking of feasible designs \\
Scale requirement & Must be preference-commensurable with the other preference scores & Need not be preference-commensurable, but should be numerically scaled for balanced search guidance \\
\bottomrule
\end{tabularx}
\caption{Comparison of acceptability relaxation and feasibility guidance.}
\end{table}

The distinction in the final two rows is essential. For negotiable limits, values beyond the original limit represent genuine actor preferences and therefore contribute to the decision-valid aggregation of alternatives. For hard constraints, guidance values merely differentiate infeasible designs during the search. Once feasibility has been reached, these values do not contribute to the preference-based ranking or recommendation.

A combined violation measure \(v\) may include violations expressed in different physical units, such as \(\text{k\texteuro}\), days, and safety units. Such a sum would not constitute a meaningful preference or decision measure, as discussed in Section~\ref{sec:limitations}. Here, however, \(v\) is used exclusively as a computational guidance measure for differentiating infeasible designs and does not enter the preference-based ranking of feasible alternatives. Its components need not be preference-commensurable, but they should be normalised or otherwise scaled to prevent one constraint from dominating the search solely because of its numerical magnitude \citep[see Axiom~2 and Example~E2]{Wolfert2026a}.

\subsection{The Need for Feasibility Guidance}

In a highly constrained problem, the feasible region may occupy only a small or fragmented part of the search space. An initial random population generated by a genetic algorithm may therefore contain few or no feasible designs. If infeasible designs are rejected immediately or all receive the same lowest score, the solver cannot distinguish candidates that are approaching the feasible region from those that remain far from it. The search then lacks a useful direction and may converge slowly, stagnate, or fail to identify a feasible design.

A feasibility-guidance function addresses this problem by assigning graded computational values to infeasible designs according to the magnitude of their constraint violations. Designs with smaller scaled violations receive a better guidance value than designs with larger violations. This allows the solver to retain and recombine candidates that move progressively towards feasibility, even when the current population contains no feasible design.

The feasibility-guidance function does not replace the hard constraint. A design is admissible only when all hard constraints are satisfied. Nor does the guidance value represent actor preference, enter the A-FINE aggregation of feasible alternatives, or justify any violation in the final recommendation. Its sole purpose is to provide computational search information on the infeasible side of the fixed constraint boundary.

The two applications may therefore use analogous mathematical constructions, but only the acceptability function represents actor preference. The feasibility-guidance function provides computational search information without decision meaning. This separation preserves the decision validity of the recommendation while helping the IMAP-IGS solver navigate highly constrained search spaces \citep{Wolfert2026b,Timp2026}. The software and code repositories are listed in the Data Availability Statement.

\subsection{Quantifying Distance to Feasibility}

Let the \(J\) hard constraints be expressed as \(h_j(\mathbf{x})\leq 0\), where \(j=1,\ldots,J\). Because the constraint violations may be expressed in different physical units, each violation is scaled before being combined. The total scaled violation of design \(\mathbf{x}\) is defined as

\begin{equation}
\label{eq:scaled-total-violation}
v(\mathbf{x})
=
\sum_{j=1}^{J}
\alpha_j
\max
\left\{
0,
\frac{h_j(\mathbf{x})}{s_j}
\right\},
\qquad
s_j>0,
\qquad
\alpha_j\geq 0,
\qquad
\sum_{j=1}^{J}\alpha_j=1,
\end{equation}

where \(s_j\) is a constraint-specific scale factor and \(\alpha_j\) controls the relative contribution of constraint \(j\) to the computational guidance measure. Consequently, \(v(\mathbf{x})=0\) if and only if design \(\mathbf{x}\) satisfies all hard constraints, whereas \(v(\mathbf{x})>0\) indicates infeasibility. The scaling factors \(s_j\) prevent constraints expressed in numerically larger units from dominating the search solely because of their magnitude.

A feasibility-guidance function \(p_v\bigl(v(\mathbf{x})\bigr)\) maps the total scaled violation to a computational guidance score between \(0\) and \(100\). It satisfies \(p_v(0)=100\) and decreases monotonically as \(v(\mathbf{x})\) increases \citep{Wolfert2026b}. Hence, all feasible designs receive the same constant guidance score, whereas infeasible designs are differentiated according to their scaled distance from feasibility. This score is used exclusively to guide the search and has no preference or decision meaning.

The feasibility-guidance function can be incorporated using constructions analogous to the two acceptability-relaxation approaches introduced in Section~\ref{sec:method}:

\begin{itemize}
    \item \textbf{Approach~1: separate feasibility-guidance criterion.} The guidance score \(p_v\bigl(v(\mathbf{x})\bigr)\) is added to the computational search score with guidance weight \(w_v\). The remaining weights are proportionally rescaled by \(1-w_v\). The resulting score is used only to compare infeasible candidate designs during the search.
    \item \textbf{Approach~2: integrated feasibility-guidance extension.} The computational scoring function associated with a constrained criterion is extended beyond the hard limit by a monotonically decreasing guidance branch \citep{Wolfert2026b}. This branch allows the solver to distinguish a small violation from a large violation, but it does not represent actor preference and cannot make an infeasible design admissible.
\end{itemize}

Although these constructions resemble the preference formulations used for negotiable acceptability limits, their interpretation is fundamentally different. For negotiable limits, scores beyond the original limit represent actor preferences and may influence the final design recommendation. For hard constraints, scores beyond the boundary are computational signals used only to guide infeasible candidates towards the feasible region.

\subsection{Constraint-dominance search rules}

Candidate designs are compared using the following constraint-dominance rules:

\begin{enumerate}
    \item Every feasible design takes precedence over every infeasible design, irrespective of their respective computational scores.
    \item Feasible designs are compared exclusively using their original A-FINE aggregate preference score \citep{Wolfert2026a, Wolfert2026b}.
    \item Infeasible designs are compared using a computational search score that combines feasibility-guidance and preference-performance information.
\end{enumerate}

Let \(\Phi(\mathbf{x})\) denote the computational search score assigned to an infeasible design. Under the separate-criterion construction of Approach~1, this score is defined as

\begin{equation}
\label{eq:infeasible-search-score}
\Phi(\mathbf{x})
=
(1-w_v)\widetilde{Z}\bigl(\mathbf{F}(\mathbf{x})\bigr)
+
w_v p_v\bigl(v(\mathbf{x})\bigr),
\qquad
0<w_v\leq 1,
\end{equation}

where \(\widetilde{Z}\bigl(\mathbf{F}(\mathbf{x})\bigr)\in[0,100]\) is the A-FINE aggregate preference score \(Z\) min-max rescaled to \([0,100]\) over the compared designs \citep[Eq.~6]{Wolfert2026a}, \(p_v\bigl(v(\mathbf{x})\bigr)\in[0,100]\) is the feasibility-guidance score, and \(w_v\) is the guidance weight. Placing both components on the same numerical interval prevents either component from dominating the search score solely because of differences in scale. The rescaled score \(\widetilde{Z}\) is used exclusively to guide the comparison of infeasible designs; the preference-based ranking of feasible designs remains based on the original A-FINE aggregate preference score \(Z\).

The feasibility-guidance term favours smaller scaled violations, while the preference-performance term retains information about the potential quality of each candidate design. A smaller violation does not necessarily take precedence over a larger violation, because preference performance and feasibility guidance jointly determine \(\Phi(\mathbf{x})\). Accordingly, \(\Phi(\mathbf{x})\) is a computational search score rather than a decision-valid aggregate preference score.

Let \(\mathcal{S}_f\) denote the feasible design space. The modified constraint-dominance relation is defined as

\begin{equation}
\label{eq:constraint-dominance}
\mathbf{x}_a \succ \mathbf{x}_b
\quad\Longleftrightarrow\quad
\begin{aligned}
&
\left[
\mathbf{x}_a\in\mathcal{S}_f
\;\land\;
\mathbf{x}_b\notin\mathcal{S}_f
\right]
\\
&{}\lor
\left[
\mathbf{x}_a,\mathbf{x}_b\in\mathcal{S}_f
\;\land\;
Z\bigl(\mathbf{F}(\mathbf{x}_a)\bigr)
>
Z\bigl(\mathbf{F}(\mathbf{x}_b)\bigr)
\right]
\\
&{}\lor
\left[
\mathbf{x}_a,\mathbf{x}_b\notin\mathcal{S}_f
\;\land\;
\Phi(\mathbf{x}_a)
>
\Phi(\mathbf{x}_b)
\right].
\end{aligned}
\end{equation}

The first condition gives every feasible design precedence over every infeasible design. The second ranks feasible designs exclusively according to their original A-FINE aggregate preference scores. The third ranks two infeasible designs according to their computational search scores, which combine feasibility-guidance and preference-performance information.

These rules constitute a modified constraint-dominance mechanism inspired by the feasibility rules of \citet{Deb2000}; see also \citet[Sec.~7.6.3]{Ehrgott2005,Martins2022}. The mechanism differs from Deb's original rules because infeasible designs are compared using the combined computational search score \(\Phi(\mathbf{x})\), rather than by total constraint violation alone. Within IMAP-IGS, this modified constraint-dominance mechanism is applied exclusively during the search and does not alter the preference-based ranking of feasible designs \citep{Timp2026}.

\subsection{Search Behaviour and Best-Fit Solution Identification}

During the initial generations of a highly constrained search, the population may contain few or no feasible designs. The feasibility-guidance score then differentiates infeasible candidates according to their scaled constraint violations and, hence, their progress towards the feasible region. Selection favours designs with better computational search scores, allowing candidates with smaller scaled violations and promising preference performance to be retained and recombined.

Once feasible designs emerge, the first constraint-dominance rule ensures that they take precedence over every infeasible design. Because \(v(\mathbf{x})=0\) for every feasible design, feasibility guidance is no longer required. Feasible designs are therefore ranked exclusively by their original A-FINE aggregate preference scores \(Z\bigl(\mathbf{F}(\mathbf{x})\bigr)\).

The final recommended design is always selected from the feasible set. If the search terminates without identifying a feasible design, the result indicates only that no feasible design was found during that search. Such an infeasible result does not constitute a design recommendation, and the hard constraints are neither relaxed nor reclassified to produce one.

The distinction between acceptability- and feasibility-constraint handling therefore lies in their meaning and use. A negotiable acceptability limit may be relaxed when justified by the resulting gain in aggregate actor preference. A hard feasibility constraint remains fixed, and its guidance function is used exclusively to differentiate infeasible candidates during the search. The final recommendation is selected exclusively from the feasible set using the original A-FINE aggregate preference score.

\section{Conclusion}

This paper has established constraint handling as an integrated preference-based design-decision problem comprising three distinct but complementary applications. First, in over-constrained stalemate problems, where no design satisfies all negotiable acceptability limits, preference-based constraint relaxation determines which limits may be relaxed and by how much to recover an acceptable solution. Second, in non-stalemate problems, where feasible designs already exist, the same approach explores whether modest departures from negotiable limits reveal substantially better designs and whether the resulting preference gain justifies those departures. In both situations, the acceptable degree of relaxation emerges endogenously from the same preference-based optimisation that identifies the best-fit design, rather than from externally imposed penalties, utility values, or violation measures.

Two equivalent formulations support these applications. Approach~1 represents each relaxation explicitly as an additional objective and search variable, providing direct visibility of the individual relaxation trade-offs. Approach~2 integrates relaxation into the preference function of each affected performance criterion, reducing the number of preference functions, weights, and search variables. When constructed consistently, both approaches yield the same best-fit solution; their practical difference lies in model complexity, interpretability, and computational effort.

Third, the same mathematical construction supports feasibility-constraint handling in highly constrained problems, but with a fundamentally different meaning. Hard constraints are never relaxed, and infeasible designs remain inadmissible. Instead, a feasibility-guidance function differentiates infeasible candidates according to their scaled constraint violations and helps the IGS solver progress towards the feasible region. Every feasible design takes precedence over every infeasible design, and feasibility guidance does not enter the final preference-based ranking. The recommended design is selected exclusively from the feasible set using the A-FINE aggregate preference score.

The resulting framework therefore unifies constraint relaxation, opportunity-seeking beyond negotiable acceptability limits, and computational guidance towards hard feasibility boundaries without conflating their distinct decision meanings. It enables preference-consistent recovery from stalemate, systematic identification of justified improvement opportunities in non-stalemate problems, and more effective IMAP-IGS search in highly constrained design spaces. In all three applications, the distinction between negotiable acceptability and non-negotiable feasibility remains explicit: negotiable limits may move only when justified by aggregate actor preference, whereas hard constraints remain fixed. The framework thereby integrates design selection, acceptability negotiation, and feasibility guidance within a single preference-based design-decision optimisation approach while preserving preference consistency and decision validity.

\section*{Data Availability Statement}

The current Preferendus software environment includes Preferendus-IMAP-IGS, which uses the A-FINE Aggregator's decision-valid preference-aggregation method to support preference-consistent design-decision optimisation. The software is publicly available through the GitHub organisations: \\ \url{https://github.com/Preferendus}\\ \url{https://github.com/TUDelft-IMAP-IGS}.

The implementations of the constraint-handling approaches presented in this paper, including the grillage-design application, are available in the Preferendus Constraint-Handling repository: \\ \url{https://github.com/Preferendus/Constraint-Handling}.

\end{document}